\documentclass[12pt]{amsart}
\usepackage{amsmath,amssymb,amsthm,mathtools,mathrsfs}
\usepackage{enumitem}
\usepackage{needspace}
\usepackage{tikz}
\usepackage[numbers,sort&compress]{natbib}
\usepackage{hyperref}
\theoremstyle{plain}
\newtheorem{theorem}{Theorem}[section]
\newtheorem{proposition}[theorem]{Proposition}
\newtheorem{lemma}[theorem]{Lemma}

\theoremstyle{remark}
\newtheorem{remark}[theorem]{Remark}
\numberwithin{equation}{section}
\allowdisplaybreaks
\newcommand{\N}{\mathbb N}
\newcommand{\R}{\mathbb R}
\newcommand{\Z}{\mathbb Z}
\newcommand{\D}{\mathscr D}
\DeclareMathOperator{\dist}{dist}
\newcommand{\dimH}{\dim_{\mathrm H}}
\newcommand{\dimB}{\dim_{\mathrm B}}
\newcommand{\dimcov}{\dim_{\mathrm{cov}}}

\makeatletter
\newcount\tm@N \newcount\tm@M \newcount\tm@halfM \newcount\tm@k \newcount\tm@t
\newcount\tm@c \newcount\tm@s \newcount\tm@pow \newcount\tm@e \newcount\tm@x
\newcount\tm@sP \newcount\tm@sA \newcount\tm@W \newcount\tm@H \newcount\tm@off
\newcount\tm@dot \newcount\tm@half
\def\tm@step{%
  \ifnum\tm@t>\z@
    \ifnum\tm@t>\tm@halfM % digit 1: c_n = c_{n-1} + s 2^{N-n+1} and T(x) = 2 - 2x
      \advance\tm@c\numexpr\tm@s*\tm@pow\relax
      \tm@s=-\tm@s
      \tm@t=\numexpr2*(\tm@M-\tm@t)\relax
    \else % digit 0: c_n = c_{n-1} and T(x) = 2x
      \tm@t=\numexpr2*\tm@t\relax
    \fi
    \tm@e=\numexpr\tm@k-\tm@c\relax % the nth error x - c_n(x), in units of 2^{-N}
    \advance\tm@sP\tm@e
    \ifnum\tm@e<\z@ \advance\tm@sA-\tm@e \else \advance\tm@sA\tm@e \fi
    \divide\tm@pow\tw@
    \expandafter\tm@step
  \fi}
\def\tm@point{%
  \tm@t=\tm@k \tm@c=\z@ \tm@s=\@ne \tm@pow=\tm@M \tm@sP=\z@ \tm@sA=\z@
  \tm@step
  \tm@x=\numexpr\tm@k*\tm@W/\tm@M-\tm@half\relax
  \pgfsys@rect{\the\tm@x sp}{\the\numexpr(\tm@sP+\tm@halfM)*\tm@H/\tm@M-\tm@half\relax sp}{\the\tm@dot sp}{\the\tm@dot sp}%
  \pgfsys@rect{\the\numexpr\tm@x+\tm@off\relax sp}{\the\numexpr\tm@sA*\tm@H/\tm@M-\tm@half\relax sp}{\the\tm@dot sp}{\the\tm@dot sp}}
\def\tm@loop{%
  \ifnum\tm@k>\tm@M\else
    \tm@point \advance\tm@k\@ne
    \expandafter\tm@loop
  \fi}
\newcommand{\tentdots}[5]{%
  \tm@N=#1\relax \tm@M=\@ne
  \tm@k=\z@ \loop\ifnum\tm@k<\tm@N \multiply\tm@M\tw@ \advance\tm@k\@ne \repeat
  \tm@halfM=\numexpr\tm@M/2\relax
  \tm@W=\dimexpr#2\relax \tm@H=\dimexpr#3\relax
  \tm@off=\dimexpr#4\relax \tm@dot=\dimexpr#5\relax \tm@half=\numexpr\tm@dot/2\relax
  \tm@k=\z@ \tm@loop \pgfsys@fill}
\makeatother

\begin{document}
\title[The signed tent-map error sum]{The signed error-sum function of tent-map expansions}
\author{Min Woong Ahn}
\address{Department of Mathematics Education, Dankook University, 152, Jukjeon-ro, Suji-gu, Yongin-si, Gyeonggi-do 16890, Republic of Korea}
\email{mwahn@dankook.ac.kr}
\date{\today}
\begin{abstract}
In this paper, we study the signed error-sum function of tent map base series introduced by Dai and Tang in 2011. We show that the function is discontinuous at every dyadic rational in the open unit interval, although it was asserted to be continuous. We determine its values at these points and analyze its range, differentiability, and variation. We also obtain sharp covering estimates for its graph and calculate the Hausdorff, box-counting, and covering dimensions of the graph, which are $1$, $1$, and $0$, respectively. Moreover, we show that the signed error sum agrees with L\'evy's function away from dyadic rationals and that the unsigned error sum is the Takagi function.
\end{abstract}
\maketitle
\tableofcontents

\section{Introduction}\label{sec:intro}
Ridley and Petruska \cite{RP00} studied an error-sum function for regular continued fractions in 2000. Error-sum functions have since been studied for several other representations, including integer-base expansions \cite{Ton07,QD10}, L\"uroth series \cite{SW07}, and Pierce expansions \cite{Ahn23}. An error-sum function is defined by summing the approximation errors of the successive convergents. If a point has two different expansions, then the value of the error sum may depend on the choice of expansion.

In 2011, Dai and Tang \cite{DT11} introduced the tent map base series, a generalized L\"uroth series in the sense of \cite{BBDK96}. To describe this expansion, let $T\colon[0,1]\to[0,1]$ be the tent map given by $T(x)\coloneqq 2\min\{x,1-x\}$. For each $n\in\N$, we denote by $T^n$ the $n$th iterate of $T$, and we put $T^0\coloneqq\operatorname{id}_{[0,1]}$. Given $x\in[0,1]$, we iterate $T$, assigning the digit $0$ to $(0,1/2]$ and the digit $1$ to $(1/2,1]$, and stop when an iterate is $0$. The resulting convergents $c_n(x)$ satisfy $|x-c_n(x)|\le2^{-n}$ (see Section~\ref{sec:prelim}), and hence the error-sum function
\begin{equation}\label{eq:P}
P(x)\coloneqq \sum_{n\in\N}(x-c_n(x))
\end{equation}
is well defined on $[0,1]$. Since the expansion of $0$ terminates at once, all convergents of $0$ are $0$, so that $P(0)=0$.

Dai and Tang asserted that $P$ is continuous \cite[Theorem~3.3 and Corollary~3.6]{DT11}. We shall show that it is discontinuous at every dyadic rational in $(0,1)$. For a dyadic rational in $(0,1)$, the two finite tent-map expansions have the same last convergent but different second-to-last convergents. This difference was not taken into account in their proof.

In this paper, we analyze the basic properties of $P$ and the dimensions of its graph, using the same digit convention as in \cite{DT11}. Put
\[
\D\coloneqq \left\{\frac{p}{2^m}:m\in\N,\ 1\le p<2^m,\ p\text{ odd}\right\},
\]
the set of dyadic rationals in $(0,1)$. The expression $r=p/2^m\in\D$ is always understood to be in lowest terms. Our first main result is concerned with the continuity, range, differentiability, and variation of $P$, and the dimensions of its graph.

\begin{theorem}\label{thm:main}
Let $P\colon[0,1]\to\R$ be the signed error-sum function defined in \eqref{eq:P}. Then the following hold.
\begin{enumerate}[label=\textup{(\roman*)},ref=\roman*,itemsep=3pt]
\item \label{thm:main(i)} The function $P$ is continuous exactly on $[0,1]\setminus\D$. For $r\coloneqq p/2^m\in\D$, the left and right limits $P(r^-)$ and $P(r^+)$ of $P$ at $r$ satisfy
\[
P(r^-)-P(r^+)=2^{1-m}\quad\text{and}\quad
P(r)=\begin{cases}
P(r^-),&\text{if }p\equiv1\pmod4,\\
P(r^+),&\text{if }p\equiv3\pmod4.
\end{cases}
\]
\item \label{thm:main(ii)} We have $P([0,1])=[-1/2,1/2]$. Moreover, $P$ has no finite derivative at any point of $(0,1)$ and has infinite total variation on every non-degenerate subinterval of $[0,1]$.
\item \label{thm:main(iii)} There exist constants $0<c<C$ such that the graph $G_P$ of $P$ satisfies
\[
c\,\delta^{-1}\log(1/\delta)\le N_\delta(G_P)\le C\,\delta^{-1}\log(1/\delta)
\]
for all sufficiently small $\delta>0$, where $N_\delta(G_P)$ denotes the least number of closed axis-parallel squares of side $\delta$ needed to cover $G_P$. In particular, $\dimH G_P=\dimB G_P=1$ and $\dimcov G_P=0$.
\end{enumerate}
\end{theorem}

We examine the argument of \cite{DT11} in Section~\ref{subsec:DT}. On the other hand, if we sum the absolute values of the errors, then the resulting function is continuous. More precisely, we have the following second main result.

\begin{theorem}\label{thm:takagi}
Let $A(x)\coloneqq \sum_{n\in\N}|x-c_n(x)|$ for $x\in[0,1]$. Then we have
\[
A(x)=\sum_{n\in\N}2^{-n}T^n(x)=\sum_{k=0}^\infty2^{-k}\dist(2^kx,\Z)\qquad(x\in[0,1]),
\]
where $\dist(y,\Z)\coloneqq\min_{j\in\Z}|y-j|$ denotes the distance from $y\in\R$ to the nearest integer. In other words, $A$ is the Takagi function. In particular, $A$ is continuous, and its graph $G_A$ satisfies
\[
\dimH G_A=\dimB G_A=\dimcov G_A=1.
\]
\end{theorem}

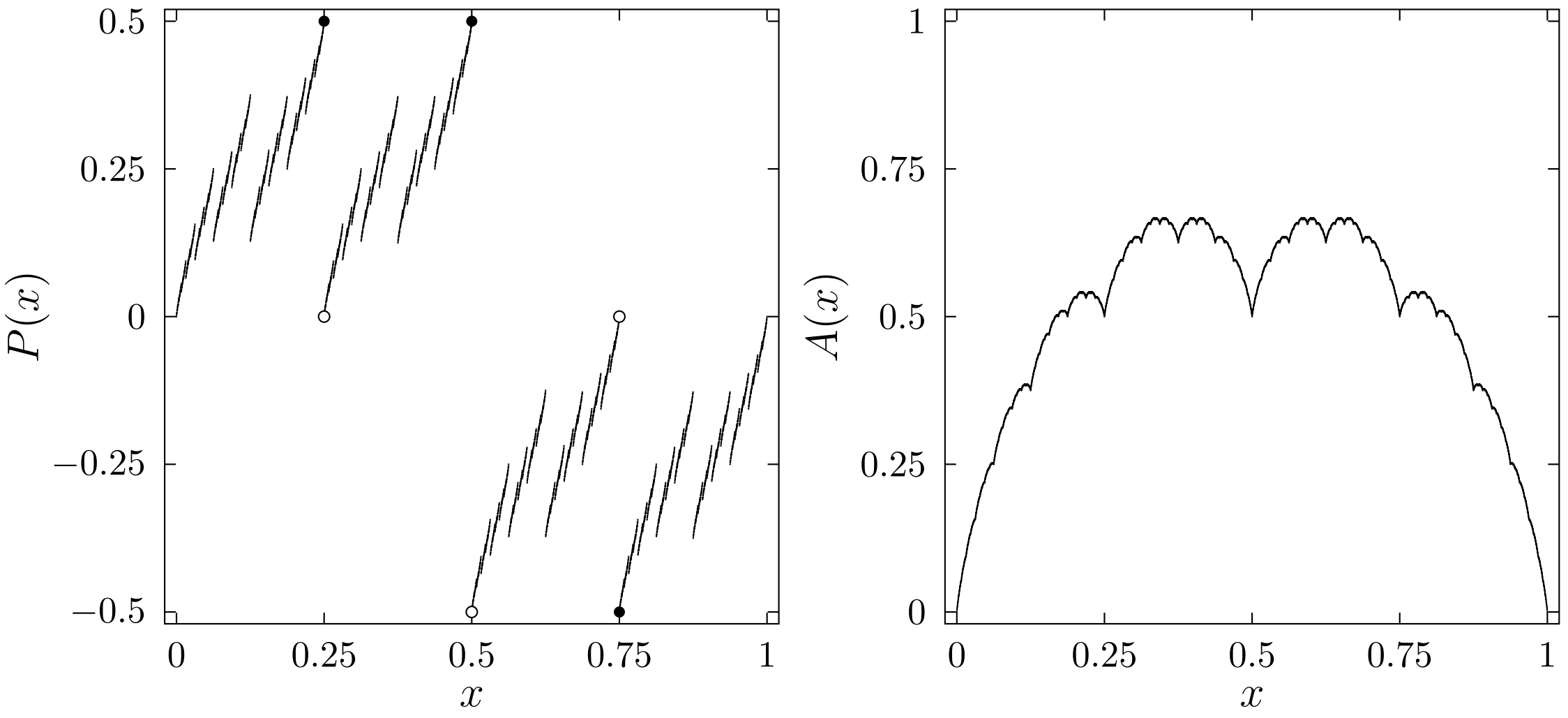
\begin{figure}[htbp]
\centering
% The points are computed and drawn by TeX itself (see \tentdots in the preamble):
% all dyadic rationals k/2^12 in [0,1], i.e., 4,097 points for each graph.
\begin{tikzpicture}[font=\footnotesize]
\def\pw{5.6cm}\def\po{7.4cm}
\foreach \ox/\lab/\yo in {0cm/{$P(x)$}/27pt, \po/{$A(x)$}/22pt} {
  \begin{scope}[xshift=\ox]
    \draw ({-0.02*\pw},{-0.02*\pw}) rectangle ({1.02*\pw},{1.02*\pw});
    \foreach \t in {0,0.25,0.5,0.75,1} {
      \draw ({\t*\pw},{-0.02*\pw}) -- ++(0,3pt) ({\t*\pw},{1.02*\pw}) -- ++(0,-3pt);
      \draw ({-0.02*\pw},{\t*\pw}) -- ++(3pt,0) ({1.02*\pw},{\t*\pw}) -- ++(-3pt,0);
      \node[anchor=north] at ({\t*\pw},{-0.02*\pw-1pt}) {$\t$};
    }
    \node[anchor=north] at ({0.5*\pw},{-0.02*\pw-13pt}) {\normalsize $x$};
    \node[anchor=south, rotate=90] at ({-0.02*\pw-\yo},{0.5*\pw}) {\normalsize \lab};
  \end{scope}
}
\foreach \t/\lb in {0/-0.5, 0.25/-0.25, 0.5/0, 0.75/0.25, 1/0.5}
  \node[anchor=east] at ({-0.02*\pw-1pt},{\t*\pw}) {$\lb$};
\foreach \t in {0,0.25,0.5,0.75,1}
  \node[anchor=east] at ({\po-0.02*\pw-1pt},{\t*\pw}) {$\t$};
\tentdots{12}{\pw}{\pw}{\po}{0.4pt}
% The value (filled) and the other one-sided limit (open) of P at x = 1/4, 1/2, 3/4.
\foreach \x/\v/\o in {0.25/1/0.5, 0.5/1/0, 0.75/0/0.5} {
  \fill ({\x*\pw},{\v*\pw}) circle (1.5pt);
  \draw[fill=white] ({\x*\pw},{\o*\pw}) circle (1.5pt);
}
\end{tikzpicture}
\caption{The graphs of $P$ (left) and of the Takagi function $A$ (right) on $[0,1]$, plotted at all dyadic rationals $k/2^{12}$ with $0\le k\le2^{12}$. The filled and open circles at $x=1/4,1/2,3/4$ indicate the value of $P$ and its other one-sided limit, respectively (see Theorem~\ref{thm:main}(\ref{thm:main(i)})).}
\label{fig:graphs}
\end{figure}

Notice that the graphs $G_P$ and $G_A$ have the same Hausdorff and box-counting dimensions, but different covering dimensions. Figure~\ref{fig:graphs} illustrates the two functions.

The function $P$ can also be expressed in terms of the binary digits. If $x\in(0,1)\setminus\D$ has the binary expansion $x=\sum_{n\in\N}b_n(x)2^{-n}$, then
\begin{equation}\label{eq:intro-binary}
P(x)=\sum_{n\in\N}(n-2)\frac{b_n(x)}{2^n}
\end{equation}
by Proposition~\ref{prop:binary}. This gives a relation between $P$ and a classical function of L\'evy. For $t\in\R$, put $\rho(t)\coloneqq t-k$ if $|t-k|<1/2$ for some $k\in\Z$, and $\rho(t)\coloneqq0$ if $t\in\Z+1/2$. Then $P$ coincides on $(0,1)\setminus\D$ with L\'evy's function
\[
\mathcal L(x)\coloneqq\sum_{n=0}^\infty\frac{\rho(2^nx)}{2^n}
\]
(see Remark~\ref{rem:base2}). Moreover, the corresponding series of absolute values agrees with the Takagi function on $(0,1)\setminus\D$, since $|\rho(t)|=\dist(t,\Z)$ for $t\notin\Z+1/2$.

L\'evy introduced $\mathcal L$ to illustrate the discontinuities of stochastic processes with independent increments (see \cite[p.~55]{KL95} and \cite{Jaf97}). Our $\mathcal L$ is the function in \cite[Section~5]{Jaf97}; the series in \cite[p.~55]{KL95} starts at $n=1$ and equals $\mathcal L(2x)/2$. In 1997, Jaffard \cite{Jaf97} described its continuity and jumps, including the jump of size $2^{1-m}$ at $p/2^m\in\D$, and determined its pointwise H\"older regularity and multifractal spectrum. More recently, Jaffard, Liao, and Zhang \cite{JLZ25,JLZ26a,JLZ26b} studied related L\'evy functions and pure jump series (see Remark~\ref{rem:base2}). It should be pointed out that $\mathcal L$ takes the average of its two one-sided limits at a dyadic rational, while $P$ takes one of them.

For the base-$2$ error sum $E_2$, we have $P(x)=E_2(x)-x$ on $(0,1)\setminus\D$ (see Remark~\ref{rem:base2}). Qiao and Dai \cite{QD10} determined the jumps of $E_2$, proved an intermediate value property, and calculated the Hausdorff dimension of its graph. Their proof of \cite[Theorem~3.1]{QD10} also gives a covering of the graph of $E_2$ by $(n+1)2^n$ squares of side $2^{-n}$. The dyadic values of $P$ and $E_2-\operatorname{id}$ differ at infinitely many points; see Remark~\ref{rem:representative}.

Thus, the assertions on the jumps and the Hausdorff dimension in Theorem~\ref{thm:main} also follow from the relation with $\mathcal L$ or $E_2$, provided that the values at dyadic rationals are taken into account. Our main purpose is to correct \cite{DT11} and study $P$ under the original digit convention. To the best of our knowledge, the relations of the tent-map error sums with L\'evy's and Takagi's functions, and the exact values at dyadic rationals, have not been recorded. We also obtain a lower bound for the covering number of the same order as the upper bound, and prove that $P$ has no finite derivative and that its graph has covering dimension $0$. These three conclusions also hold for $\mathcal L$ (Remark~\ref{rem:representative}). The assertion on the covering dimension can also be proved by using the binary representation, which shows that the closure of $G_P$ is homeomorphic to the Cantor space $\{0,1\}^\N$ (Remark~\ref{rem:clopen}). All parts of Theorem~\ref{thm:main} are proved directly by elementary arguments.

This paper is organized as follows. In Section~\ref{sec:prelim}, we recall the tent map base series and establish some elementary formulas for $P$. In Section~\ref{sec:aux}, we prove three auxiliary lemmas needed in the proofs of the main results. In Section~\ref{sec:proofs}, we prove Theorems~\ref{thm:main} and~\ref{thm:takagi} and examine the argument of \cite{DT11}. Finally, Section~\ref{sec:remarks} contains further remarks on the choice of dyadic values and some alternative proofs.

Throughout the paper, we denote by $\N$ the set of positive integers. All subsets of Euclidean spaces are endowed with the usual subspace topology. For a function $f$ on $[0,1]$, we write $f(x^-)$ and $f(x^+)$ for its left and right limits at $x$, whenever they exist. We say that $f$ is {\it regulated} if $f(x^-)$ exists for every $x\in(0,1]$ and $f(x^+)$ exists for every $x\in[0,1)$. The graph of $f$ is denoted by $G_f\coloneqq\{(x,f(x)):x\in[0,1]\}$.

For a bounded set $F\subseteq\R^2$, let $N_\delta(F)$ denote the least number of closed axis-parallel squares of side $\delta>0$ needed to cover $F$. We denote by $\dimH F$ the Hausdorff dimension, by $\dimB F$ the box-counting dimension when it exists, and by $\dimcov F$ the Lebesgue covering dimension. We refer the reader to \cite[Chapters~2 and~3]{Fal14} and \cite{Eng95} for details on these dimensions. The logarithm $\log$ always means the natural logarithm.

\section{Preliminaries}\label{sec:prelim}
In this section, we recall the definition of the tent map base series from \cite{DT11} and present some elementary facts about its convergents and error sum.

For $x\in[0,1]$, put $x_0\coloneqq x$ and $x_n\coloneqq T^n(x)$ for $n\in\N$. As long as $x_{n-1}\ne0$, define
\[
a_n(x)\coloneqq \begin{cases}0,&\text{if }0<x_{n-1}\le1/2,\\1,&\text{if }1/2<x_{n-1}\le1.\end{cases}
\]
For $x>0$, let $N(x)$ be the least $n\in\N$ such that $x_n=0$, with $N(x)\coloneqq\infty$ if no such $n$ exists. For each $n\in\N$ with $n\le N(x)$, define
\begin{equation}\label{eq:convergents}
c_n(x)\coloneqq \sum_{j=1}^n(-1)^{a_1(x)+\cdots+a_{j-1}(x)}\frac{2a_j(x)}{2^j},
\end{equation}
where the exponent is $0$ for $j=1$, and put $c_n(x)\coloneqq x$ for $n>N(x)$. Moreover, put $c_n(0)\coloneqq 0$ for every $n\in\N$. We call $c_n(x)$ the $n$th convergent of $x$, and we write $e_n(x)\coloneqq x-c_n(x)$ for the $n$th error.

For $x>0$, we have $x_{j-1}=a_j+(-1)^{a_j}x_j/2$ for $j\le N(x)$. Iterating this relation, we obtain
\begin{equation}\label{eq:remainder}
e_n(x)=(-1)^{a_1(x)+\cdots+a_n(x)}2^{-n}T^n(x)
\end{equation}
whenever $a_1(x),\dots,a_n(x)$ are defined. When the expansion terminates, both $e_n(x)$ and $T^n(x)$ are $0$ from that point on. Thus, for every $x\in[0,1]$ and $n\in\N$, we have $|e_n(x)|=2^{-n}T^n(x)\le2^{-n}$. It follows that the series defining $P$ converges absolutely and uniformly on $[0,1]$. In particular, the partial sums $S_N(x)\coloneqq\sum_{n=1}^Ne_n(x)$ satisfy
\begin{equation}\label{eq:tail}
\sup_{x\in[0,1]}|P(x)-S_N(x)|\le2^{-N}\qquad(N\in\N).
\end{equation}

For later use, we give an estimate on dyadic intervals. For $n\in\N$ and $0\le k<2^n$, put $I_{n,k}\coloneqq(k2^{-n},(k+1)2^{-n})$. If $x\in I_{n,k}$, then $T^j(x)\notin\{0,1/2\}$ for $0\le j\le n-1$. By the continuity of $T^j$ and the connectedness of $I_{n,k}$, each $T^j(I_{n,k})$ lies entirely in one of the two digit intervals. Hence, the digits $a_1,\dots,a_n$, and consequently the convergents $c_1,\dots,c_n$, are constant on $I_{n,k}$. It follows that $S_n$ is affine with slope $n$ there. Combining this with \eqref{eq:tail}, we obtain
\begin{equation}\label{eq:cell}
|P(y)-P(x)-n(y-x)|\le2^{1-n}\qquad(x,y\in I_{n,k}).
\end{equation}
We will use \eqref{eq:cell} to prove the assertion on the variation and the covering estimates in Theorem~\ref{thm:main}. Moreover, each $c_n$ is a step function, so that each $e_n$ is regulated. By uniform convergence, $P$ is also regulated, and its one-sided limits may be computed term by term.

The following proposition gives a formula for $P$ in terms of the ordinary binary digits.

\begin{proposition}\label{prop:binary}
Let $x\in(0,1)\setminus\D$ have the binary expansion $x=\sum_{j\in\N}b_j2^{-j}$. Then we have
\begin{equation}\label{eq:binary-conv}
c_n(x)=\sum_{j=1}^{n-1}\frac{b_j}{2^j}+\frac{2b_n}{2^n}\quad(n\in\N)
\qquad\text{and}\qquad
P(x)=\sum_{j\in\N}(j-2)\frac{b_j}{2^j}.
\end{equation}
\end{proposition}

\begin{proof}
Put $\sigma_k\coloneqq a_1(x)+\cdots+a_k(x)$ for $k\ge0$, with $\sigma_0\coloneqq 0$, and put $b_0\coloneqq 0$. Since $x\notin\D$, the expansion does not terminate and $T^k(x)\ne1/2$ for every $k\ge0$. For $y\in(0,1)\setminus\D$ with binary digits $y_1,y_2,\dots$, we have $T(y)=2y=\sum_{i\in\N}y_{i+1}2^{-i}$ if $y_1=0$, and $T(y)=2-2y=\sum_{i\in\N}(1-y_{i+1})2^{-i}$ if $y_1=1$, since $\sum_{i\in\N}2^{-i}=1$. In other words, $T$ deletes the first binary digit and, if this digit is $1$, also complements the remaining digits. Using this observation, we find by induction on $k$ that the binary digits of $T^k(x)$ are
\[
b_{k+i}+\sigma_k\pmod2\qquad(i\in\N).
\]
Comparing the first digits, we obtain $a_{k+1}(x)\equiv b_{k+1}+\sigma_k\pmod2$. Thus, $\sigma_{k+1}\equiv b_{k+1}\pmod2$, so that $\sigma_k\equiv b_k\pmod2$ for all $k\ge0$. Therefore, $a_j(x)\equiv b_{j-1}+b_j\pmod2$, i.e., $a_j(x)=|b_j-b_{j-1}|$, and since $(-1)^{\sigma_{j-1}}=(-1)^{b_{j-1}}$, we obtain
\[
(-1)^{\sigma_{j-1}}a_j(x)=b_j-b_{j-1}\qquad(j\in\N).
\]
Substituting this into \eqref{eq:convergents}, we obtain $c_n(x)=\sum_{j=1}^n(b_j-b_{j-1})2^{1-j}$, and collecting the terms containing each $b_j$ gives the first formula in \eqref{eq:binary-conv}. Subtracting it from $x=\sum_{j\in\N}b_j2^{-j}$, we obtain
\[
e_n(x)=\sum_{j=n+1}^\infty\frac{b_j}{2^j}-\frac{b_n}{2^n}.
\]
Summing over $n$, we obtain
\[
P(x)=\sum_{n\in\N}\sum_{j=n+1}^\infty\frac{b_j}{2^j}-\sum_{n\in\N}\frac{b_n}{2^n}.
\]
In the double series, all terms are non-negative, and each $b_j2^{-j}$ appears for exactly $j-1$ values of $n$, namely $1\le n\le j-1$. Since $\sum_{j\in\N}j2^{-j}<\infty$, changing the order of summation gives
\[
\sum_{n\in\N}\sum_{j=n+1}^\infty\frac{b_j}{2^j}=\sum_{j\in\N}(j-1)\frac{b_j}{2^j}.
\]
The second sum is $x=\sum_{j\in\N}b_j2^{-j}$, and so we obtain the second formula in \eqref{eq:binary-conv}.
\end{proof}

\begin{remark}\label{rem:base2}
For $x\in[0,1)$, let $x=\sum_{j\in\N}b_j(x)2^{-j}$ be the binary expansion of $x$, where the expansion of a dyadic rational is taken to be terminating, and let $d_n(x)\coloneqq \sum_{j=1}^nb_j(x)2^{-j}$ be the ordinary binary convergents. Let $E_2(x)\coloneqq \sum_{n\in\N}(x-d_n(x))$ be the base-$2$ error sum. Since $x-d_n(x)=2^{-n}\{2^nx\}$, where $\{\cdot\}$ denotes the fractional part, changing the order of summation as in the proof of Proposition~\ref{prop:binary} gives
\[
E_2(x)=\sum_{n\in\N}\frac{\{2^nx\}}{2^n}=\sum_{j\in\N}(j-1)\frac{b_j(x)}{2^j}.
\]
Moreover, Proposition~\ref{prop:binary} gives
\[
P(x)=E_2(x)-x\qquad(x\in(0,1)\setminus\D).
\]
We next verify the relation with L\'evy's function. Let $\rho$ and $\mathcal L$ be as in Section~\ref{sec:intro}, and let $x\in(0,1)\setminus\D$. For every integer $n\ge0$, we have $\{2^nx\}=\sum_{i\in\N}b_{n+i}(x)2^{-i}$, and since $x\notin\D$, we have $\{2^nx\}\ne1/2$. Moreover, $\{2^nx\}<1/2$ if and only if $b_{n+1}(x)=0$. Hence, the integer nearest to $2^nx$ is $\lfloor2^nx\rfloor+b_{n+1}(x)$, so that
\[
\rho(2^nx)=\{2^nx\}-b_{n+1}(x)\qquad(n\ge0).
\]
Since $\sum_{n=0}^\infty2^{-n}\{2^nx\}=x+E_2(x)$ and $\sum_{n=0}^\infty2^{-n}b_{n+1}(x)=2x$, we obtain
\[
\mathcal L(x)=\sum_{n=0}^\infty\frac{\{2^nx\}-b_{n+1}(x)}{2^n}=\bigl(x+E_2(x)\bigr)-2x=P(x).
\]
At each point of $\D$, however, $\mathcal L$ takes the average of its two one-sided limits, since $\rho$ vanishes at the half-integers, whereas $P$ takes one of them. We also note that, on $(0,1)\setminus\D$, the function $E_2-1/2$ coincides with the lacunary Davenport series $\sum_{n\in\N}2^{-n}(\{2^nx\}-1/2)$ (see \cite{Jaf04} and \cite[Definition~4.4]{JLZ25}) and with the L\'evy function $L^2_1$ of \cite{JLZ26a}, so that $L^2_1(x)=\mathcal L(x)+x-1/2$ for $x\in(0,1)\setminus\D$. For further results on the pointwise regularity and multifractal analysis of pure jump series, including this lacunary Davenport series, see \cite{JLZ25}. We point out that these identities hold away from $\D$ and do not determine the values of $P$ at dyadic rationals. For these values, we use the tent-map expansion defined above.
\end{remark}

\section{Auxiliary results}\label{sec:aux}
In this section, we establish three lemmas needed in Section~\ref{sec:proofs}. We begin with the behavior of the individual errors near a dyadic rational.

\begin{lemma}\label{lem:errors}
Let $r\coloneqq p/2^m\in\D$. Then the following hold.
\begin{enumerate}[label=\textup{(\roman*)},ref=\roman*,itemsep=2pt]
\item \label{lem:errors(i)} We have $T^{m-1}(r)=1/2$, $T^m(r)=1$, and $T^{m+1}(r)=0$. Moreover, $a_m(r)=0$, and the digits $a_1,\dots,a_{m-1}$ are constant on some open interval containing $r$.
\item \label{lem:errors(ii)} For $n<m$, the function $e_n$ is continuous at $r$.
\item \label{lem:errors(iii)} For $n>m$, we have $e_n(r^-)=e_n(r)=e_n(r^+)=0$.
\item \label{lem:errors(iv)} We have $e_m(r^-)=2^{-m}$, $e_m(r^+)=-2^{-m}$, and $e_m(r)=(-1)^{(p-1)/2}2^{-m}$.
\end{enumerate}
\end{lemma}

\begin{proof}
(\ref{lem:errors(i)}) Note that $T$ maps a reduced fraction $q/2^l$ with $l\in\N$ to a reduced fraction with denominator $2^{l-1}$. Hence, $T^{m-1}(r)=1/2$, $T^m(r)=1$, and $T^{m+1}(r)=0$, with $T^j(r)\notin\{0,1/2,1\}$ for $0\le j\le m-2$. Consequently, $a_m(r)=0$. Moreover, for $0\le j\le m-2$, the point $T^j(r)$ lies in $(0,1/2)$ or in $(1/2,1)$, and so, by the continuity of $T^j$, so does $T^j(x)$ for all $x$ in some open interval containing $r$. Since the digit $a_{j+1}(x)$ is determined by the position of $T^j(x)$, the digits $a_1,\dots,a_{m-1}$ are constant on some open interval containing $r$.

(\ref{lem:errors(ii)}) Let $n<m$, and let $U$ be an open interval containing $r$ on which $a_1,\dots,a_{m-1}$ are constant, as in part~(\ref{lem:errors(i)}). Since $n\le m-1$, the sign $\varepsilon\coloneqq(-1)^{a_1(x)+\cdots+a_n(x)}$ does not depend on $x\in U$, and so \eqref{eq:remainder} gives $e_n(x)=\varepsilon2^{-n}T^n(x)$ for $x\in U$. Since $T^n$ is continuous, $e_n$ is continuous at $r$.

(\ref{lem:errors(iii)}) Let $n>m$. By part~(\ref{lem:errors(i)}), we have $T^n(r)=0$, and so $e_n(r)=0$ since $|e_n|=2^{-n}T^n$. Moreover, since $T^n$ is continuous, $|e_n(x)|=2^{-n}T^n(x)$ tends to $0$ as $x\to r$, which gives $e_n(r^-)=e_n(r^+)=0$.

(\ref{lem:errors(iv)}) Since $p$ is odd, the terminating binary expansion of $r$ has exactly $m$ digits, i.e., $r=\sum_{j=1}^mb_j2^{-j}$ with $b_m=1$. In particular, $\sum_{j=1}^{m-1}b_j2^{-j}=r-2^{-m}$. The convergent $c_m$ is constant on each of the intervals $(r-2^{-m},r)$ and $(r,r+2^{-m})$. Applying Proposition~\ref{prop:binary} at points of $(0,1)\setminus\D$ in these intervals, we find that $c_m=r-2^{-m}$ on the first interval and $c_m=r+2^{-m}$ on the second. Hence, $e_m(r^-)=2^{-m}$ and $e_m(r^+)=-2^{-m}$. On the other hand, by part~(\ref{lem:errors(i)}) and \eqref{eq:remainder}, we have
\[
e_m(r)=(-1)^{a_1(r)+\cdots+a_{m-1}(r)}2^{-m}.
\]
On the open interval in part~(\ref{lem:errors(i)}), the map $T^{m-1}$ is a composition of the branches $y\mapsto2y$ and $y\mapsto2-2y$, the latter used exactly when the corresponding digit is $1$. Hence, $T^{m-1}$ has slope $(-1)^{a_1(r)+\cdots+a_{m-1}(r)}2^{m-1}$ near $r$. On the other hand, the map $T^{m-1}$ is affine with slope $\pm2^{m-1}$ on each interval $[k2^{1-m},(k+1)2^{1-m}]$, and it maps the points $k2^{1-m}$ alternately to $0$ and $1$, starting with $T^{m-1}(0)=0$. Thus, the sign of the slope on this interval is $(-1)^k$. Since $p$ is odd, $r$ is the midpoint of this interval for $k\coloneqq (p-1)/2$. Therefore, $(-1)^{a_1(r)+\cdots+a_{m-1}(r)}=(-1)^{(p-1)/2}$, and so $e_m(r)=(-1)^{(p-1)/2}2^{-m}$.
\end{proof}

The following lemma will be used to determine the range of $P$. It is a variant of \cite[Theorem~3.23]{Ahn23}, in that the value at a discontinuity is allowed to lie strictly between the two one-sided limits. We include an elementary proof for the convenience of the reader.

\begin{lemma}\label{lem:ivp}
Let $f\colon[0,1]\to\R$ be a regulated function satisfying
\begin{equation}\label{eq:ordered-f}
f(x^-)\ge f(x)\ge f(x^+)\qquad(x\in(0,1)).
\end{equation}
If $0<a<b<1$ and $f(a)<y<f(b)$, then there exists $c\in(a,b)$ such that $f$ is continuous at $c$ and $f(c)=y$. The same conclusion holds when $a=0$ or $b=1$, provided $f$ is continuous at each endpoint used.
\end{lemma}

\begin{proof}
Put $c\coloneqq \sup\{x\in[a,b]:f(x)<y\}$. By \eqref{eq:ordered-f} and the continuity assumption when $a=0$ or $b=1$, we have $f(a^+)\le f(a)<y$ and $f(b^-)\ge f(b)>y$. Hence, $a<c<b$. If $f(c)<y$, then $f(c^+)\le f(c)<y$, so that $f(x)<y$ for some $x\in(c,b)$. This contradicts the definition of $c$. Hence, $f(c)\ge y$, and there are points of the set defining $c$ which tend to $c$ from the left. Consequently, $f(c^-)\le y$. On the other hand, $f(c^+)\ge y$ by the definition of $c$. Combining these inequalities with \eqref{eq:ordered-f}, we obtain $f(c^-)=f(c)=f(c^+)=y$. Hence, $f$ is continuous at $c$ and $f(c)=y$, as required.
\end{proof}

Finally, recall that the covering dimension $\dimcov X$ of a topological space $X$ is the least integer $n\ge-1$ such that every finite open cover of $X$ has a finite open refinement in which every point of $X$ lies in at most $n+1$ members (see \cite{Eng95}). In particular, a non-empty space $X$ satisfies $\dimcov X=0$ if and only if every finite open cover of $X$ has a finite refinement consisting of pairwise disjoint open sets. We shall use the following lemma to prove that $\dimcov G_P=0$. For this purpose, we will construct suitable finite partitions of $G_P$.

\begin{lemma}\label{lem:partitions}
Let $X$ be a non-empty metric space. For each $n\in\N$, let $\mathscr P_n$ be a finite partition of $X$ into open sets. Suppose that $\mathscr P_{n+1}$ refines $\mathscr P_n$ for every $n\in\N$, and that the largest diameter of the members of $\mathscr P_n$ tends to $0$ as $n\to\infty$. Then $\dimcov X=0$.
\end{lemma}

\begin{proof}
Put $\mathscr P\coloneqq\bigcup_{n\in\N}\mathscr P_n$. We first note two facts which follow from the refinement assumption. Any two members of $\mathscr P$ are either disjoint or one contains the other. Indeed, if $A\in\mathscr P_i$ and $B\in\mathscr P_j$ with $i\le j$, then $B$ is contained in a member of $\mathscr P_i$, which either equals $A$ or is disjoint from $A$. Moreover, if $B\in\mathscr P_n$, then every member of $\mathscr P_i$ with $i\ge n$ that contains $B$ is equal to $B$. For each $i<n$, there is only one member of $\mathscr P_i$ containing $B$. Hence, $B$ is contained in only finitely many distinct members of $\mathscr P$.

Let $\mathscr U\coloneqq\{U_1,\dots,U_\ell\}$ be a finite open cover of $X$. For convenience, we call a member of $\mathscr P$ {\it small} if it is contained in some $U_i$. Let $x\in X$, and choose $i$ and $\rho>0$ such that the open ball of radius $\rho$ centered at $x$ is contained in $U_i$. For all sufficiently large $n$, the member of $\mathscr P_n$ containing $x$ has diameter less than $\rho$, and so it is small. Thus, the small sets cover $X$. By the previous paragraph, every small set is contained in a maximal small set, and two distinct maximal small sets are disjoint. For $1\le i\le\ell$, let $V_i$ be the union of those maximal small sets that are contained in $U_i$ but in none of $U_1,\dots,U_{i-1}$. Then $V_1,\dots,V_\ell$ are pairwise disjoint open sets covering $X$, and $V_i\subseteq U_i$ for each $i$. Hence, $\dimcov X\le0$. Since $X$ is non-empty, we conclude that $\dimcov X=0$.
\end{proof}

\section{Proofs of main results}\label{sec:proofs}
In this section, we prove the main results stated in Section~\ref{sec:intro}. We also examine the argument of \cite{DT11} in Section~\ref{subsec:DT}.

\subsection{Proof of Theorem~\ref{thm:main}(\ref{thm:main(i)})}\label{subsec:i}
\begin{proof}[Proof of Theorem~\ref{thm:main}(\ref{thm:main(i)})]
We first show that $P$ is continuous on $[0,1]\setminus\D$. Let $x\in(0,1)\setminus\D$. For each $n\in\N$, the point $x$ belongs to some $I_{n,k}$, on which $c_n$ is constant. Hence, each $e_n$ is continuous at $x$. By the uniform convergence in \eqref{eq:tail}, $P$ is continuous at $x$ as well. For the endpoints, note that $T^n(0)=T^n(1)=0$ for $n\in\N$. Since $|e_n|=2^{-n}T^n$, each $e_n$ tends to $0$ as $x\to0^+$ or $x\to1^-$. Again, by uniform convergence, $P$ is continuous at $0$ and $1$, with $P(0)=P(1)=0$.

Now, let $r\coloneqq p/2^m\in\D$. By Lemma~\ref{lem:errors}(\ref{lem:errors(ii)}) and~(\ref{lem:errors(iii)}), the values $e_n(r^-)$, $e_n(r)$, and $e_n(r^+)$ agree for every $n\ne m$. Thus, when we calculate the differences among $P(r^-)$, $P(r)$, and $P(r^+)$ term by term, all terms except the $m$th cancel. By Lemma~\ref{lem:errors}(\ref{lem:errors(iv)}), we obtain
\[
P(r^-)-P(r)=2^{-m}-(-1)^{(p-1)/2}2^{-m}
\quad\text{and}\quad
P(r)-P(r^+)=(-1)^{(p-1)/2}2^{-m}+2^{-m}.
\]
Hence, $P(r^-)-P(r^+)=2^{1-m}$. Moreover, $P(r)=P(r^-)$ if $p\equiv1\pmod4$, and $P(r)=P(r^+)$ if $p\equiv3\pmod4$. In particular, $P$ is discontinuous at $r$.
\end{proof}

For later use, we record explicit formulas for the one-sided limits. Let $r\coloneqq p/2^m\in\D$, and let $r=\sum_{j=1}^mb_j2^{-j}$ be its terminating binary expansion. As $x\in(0,1)\setminus\D$ tends to $r$ from the right, each binary digit $b_j(x)$ is eventually equal to $b_j$ if $1\le j\le m$, and to $0$ if $j>m$. Since $\sum_{j\in\N}|j-2|2^{-j}<\infty$, Proposition~\ref{prop:binary} and Theorem~\ref{thm:main}(\ref{thm:main(i)}) give
\begin{equation}\label{eq:traces}
P(r^+)=\sum_{j=1}^m(j-2)\frac{b_j}{2^j}
\quad\text{and}\quad
P(r^-)=P(r^+)+2^{1-m}.
\end{equation}
In particular, \eqref{eq:traces} and Theorem~\ref{thm:main}(\ref{thm:main(i)}) give
\begin{equation}\label{eq:special}
P(1/4)=\frac12,\qquad P(1/2)=\frac12,\qquad P(3/4)=-\frac12.
\end{equation}
These values can also be checked directly from the definition. In fact, the convergents of $1/4$ are $0,0,1/4,\dots$, those of $1/2$ are $0,1/2,\dots$, and those of $3/4$ are $1,1,3/4,\dots$.

\subsection{The argument in \cite{DT11}}\label{subsec:DT}
We first calculate $P$ near $1/2$ directly from the convergents. For $k\ge3$, put $y_k^\pm\coloneqq 1/2\pm2^{-k}$. The first convergents of $y_k^-$ and $y_k^+$ are $0$ and $1$, respectively. Their convergents $c_n$ with $2\le n\le k$ are equal to $1/2$, and their convergents $c_n$ with $n\ge k+1$ are equal to the points themselves. Hence, we have
\[
P(y_k^-)=\frac12-k2^{-k}
\quad\text{and}\quad
P(y_k^+)=-\frac12+k2^{-k}.
\]
Since $P$ is regulated, it follows that $P((1/2)^-)=1/2$ and $P((1/2)^+)=-1/2$. In particular, changing the value of $P$ at $1/2$ cannot make it continuous there.

We now compare the two finite expansions used in \cite[Theorem~3.3, pp.~575--576]{DT11}. Let $n\ge2$, and let $u\coloneqq(u_1,\dots,u_{n-2})\in\{0,1\}^{n-2}$ be a word of length $n-2$, which is empty if $n=2$. Put
\[
C\coloneqq\sum_{j=1}^{n-2}(-1)^{u_1+\cdots+u_{j-1}}\frac{2u_j}{2^j}
\quad\text{and}\quad
s\coloneqq(-1)^{u_1+\cdots+u_{n-2}},
\]
so that $C=0$ and $s=1$ if $n=2$. The two words $(u,0,1)$ and $(u,1,1)$ represent the same point. Substituting each word into \eqref{eq:convergents}, we obtain the following last two convergents:
\[
\begin{array}{c|cc}
 & c_{n-1}&c_n\\ \hline
(u,0,1)&C&C+s2^{1-n}\\
(u,1,1)&C+s2^{2-n}&C+s2^{1-n}
\end{array}
\]
All corresponding convergents agree except the second-to-last ones, which differ by $s2^{2-n}$. Hence, the two finite error sums differ by $s2^{2-n}$.

In Case~I of the proof of \cite[Theorem~3.3]{DT11}, we have $s=1$. The endpoint denoted there by $x_1$, which we write as $r$, uses the word $(u,0,1)$, while the points $r+\alpha$ considered there use the prefix $(u,1,1)$. In the displayed calculation for $P(r+\alpha)$ on p.~576, the first $n$ convergents of $r+\alpha$ are replaced by those of $r$. This replacement is incorrect for $c_{n-1}$. Indeed, the first $n$ errors of $r+\alpha$ sum to
\[
P(r)+n\alpha-2^{2-n},
\]
not to $P(r)+n\alpha$. The omitted term $2^{2-n}$ is precisely the jump of $P$ at $r$ given by Theorem~\ref{thm:main}(\ref{thm:main(i)}). Therefore, \cite[Theorem~3.3 and Corollary~3.6]{DT11} are false for the signed sum defined there.

The bound in \cite[Proposition~3.4]{DT11} is also false for the values given by the original digit convention. In the notation of \cite{DT11}, the word $(1)$ gives the cylinder $(1/2,1]$, whose left endpoint is $1/2$. Thus, the lower bound in \cite[Proposition~3.4]{DT11}, applied with $x=3/4$, would give $P(1/2)<P(3/4)$, which contradicts \eqref{eq:special}. The continuity of $P$ on $[0,1]\setminus\D$ follows directly from uniform convergence, as shown in Section~\ref{subsec:i}.

\subsection{Proof of Theorem~\ref{thm:main}(\ref{thm:main(ii)})}\label{subsec:ii}
\begin{proof}[Proof of Theorem~\ref{thm:main}(\ref{thm:main(ii)})]
By Theorem~\ref{thm:main}(\ref{thm:main(i)}), the function $P$ is regulated, continuous at $0$ and $1$, and satisfies
\[
P(x^-)\ge P(x)\ge P(x^+)\qquad(x\in(0,1)).
\]
Thus, $P$ satisfies the hypotheses of Lemma~\ref{lem:ivp}.

We first determine its range. For $x\in(0,1)\setminus\D$, Proposition~\ref{prop:binary} gives $P(x)\ge-b_1/2\ge-1/2$, since only the term with $j=1$ is negative, and
\[
P(x)\le\sum_{j=3}^\infty(j-2)2^{-j}=\frac12.
\]
At each point of $\D$, the value of $P$ is one of its one-sided limits by Theorem~\ref{thm:main}(\ref{thm:main(i)}). Together with $P(0)=P(1)=0$, this shows that the same bound holds on $[0,1]$. By \eqref{eq:special}, the values $-1/2$ and $1/2$ are both attained. Applying Lemma~\ref{lem:ivp} on $[0,1/4]$, we see that every value in $(0,1/2)$ is attained. Similarly, applying Lemma~\ref{lem:ivp} on $[3/4,1]$, we see that every value in $(-1/2,0)$ is attained. Since $P(0)=0$, we conclude that $P([0,1])=[-1/2,1/2]$.

We next prove that the total variation is infinite on every non-degenerate interval. Fix a closed interval $[a,b]\subseteq[0,1]$ with $a<b$, and let $n\ge5$. For each $k$ with $[k2^{-n},(k+1)2^{-n}]\subseteq[a,b]$, put $u_k\coloneqq(k+1/4)2^{-n}$ and $v_k\coloneqq(k+3/4)2^{-n}$. Then $u_k,v_k\in I_{n,k}$ and $v_k-u_k=2^{-n-1}$, and so \eqref{eq:cell} gives
\[
P(v_k)-P(u_k)\ge n2^{-n-1}-2^{1-n}=(n-4)2^{-n-1}.
\]
There are at least $(b-a)2^n-2$ such $k$. Take a partition of $[a,b]$ containing all these points $u_k$ and $v_k$ in increasing order. The variation sum for this partition is at least $((b-a)2^n-2)(n-4)2^{-n-1}$, which tends to infinity as $n\to\infty$. Hence, $P$ has infinite total variation on $[a,b]$. Since every non-degenerate subinterval of $[0,1]$ contains such a closed interval, the assertion on the variation follows. In particular, since $P$ is bounded, $P$ is neither non-decreasing nor non-increasing on any non-degenerate interval.

Lastly, we show that $P$ has no finite derivative. For a point in $\D$, this follows from the discontinuity of $P$. Let $x\in(0,1)\setminus\D$ have binary digits $(b_j)_{j\in\N}$, and let $y_n$ be the point obtained by changing only the $n$th binary digit of $x$. Then $y_n\in(0,1)\setminus\D$, $y_n\to x$ as $n\to\infty$, and Proposition~\ref{prop:binary} gives
\[
y_n-x=(1-2b_n)2^{-n}
\quad\text{and}\quad
\frac{P(y_n)-P(x)}{y_n-x}=n-2.
\]
Since the last quotient tends to infinity as $n\to\infty$, $P$ has no finite derivative at $x$.
\end{proof}

\subsection{Proof of Theorem~\ref{thm:main}(\ref{thm:main(iii)})}\label{subsec:iii}
\begin{proof}[Proof of Theorem~\ref{thm:main}(\ref{thm:main(iii)})]
We first prove the upper and lower bounds for the covering number. For the upper bound, let $n\in\N$ and put $\delta\coloneqq2^{-n}$. By \eqref{eq:cell}, the values of $P$ on each $I_{n,k}$, $0\le k<2^n$, lie in a closed interval of length $(n+2)\delta$. Thus, the part of $G_P$ above $I_{n,k}$ lies in a rectangle of width $\delta$ and height $(n+2)\delta$, so it can be covered by $n+2$ squares of side $\delta$. Covering the $2^n+1$ points of $G_P$ above the endpoints $k2^{-n}$ separately, we obtain
\begin{equation}\label{eq:upper}
N_{2^{-n}}(G_P)\le(n+2)2^n+2^n+1\le4(n+1)2^n.
\end{equation}
Notice that we have covered the dyadic endpoints separately. Hence, \eqref{eq:upper} holds regardless of which one-sided limit is chosen as the value at each endpoint.

For the lower bound, let $n\ge4$, put $\delta\coloneqq2^{-n}$ and $L\coloneqq\lfloor n/4\rfloor$, and fix $0\le k<2^n$. For $1\le i\le L$, put $t_i\coloneqq(k+(4i-2)/n)\delta$. Since $4L-2<n$, these points lie in $I_{n,k}$, and \eqref{eq:cell} gives
\[
P(t_j)-P(t_i)\ge4(j-i)\delta-2\delta\ge2\delta\qquad(1\le i<j\le L).
\]
Hence, no closed square of side $\delta$ contains two of the points $(t_i,P(t_i))$, $1\le i\le L$. Moreover, the horizontal projection of a closed square of side $\delta$ meets at most two of the intervals $I_{n,k}$, $0\le k<2^n$. Choosing $L$ points in this way for each of the $2^n$ intervals $I_{n,k}$, we obtain $L2^n$ points of $G_P$. By the preceding two observations, a closed square of side $\delta$ contains at most two of these points. Since $L\ge n/8$, it follows that
\[
N_{2^{-n}}(G_P)\ge\frac{L2^n}{2}\ge\frac{n2^n}{16}.
\]
Together with \eqref{eq:upper}, this proves the covering estimate for $\delta=2^{-n}$. For $0<\delta\le2^{-4}$, choose $n\ge4$ such that $2^{-(n+1)}<\delta\le2^{-n}$. Since $N_\delta(G_P)$ is non-increasing in $\delta$, we have $N_{2^{-n}}(G_P)\le N_\delta(G_P)\le N_{2^{-(n+1)}}(G_P)$. Moreover, $(\log2)\,n2^n\le\delta^{-1}\log(1/\delta)<(4\log2)\,n2^n$. Hence, the bounds just proved give
\[
\frac{1}{64\log2}\,\delta^{-1}\log(1/\delta)\le N_\delta(G_P)\le\frac{24}{\log2}\,\delta^{-1}\log(1/\delta).
\]

Consequently, the lower and upper box-counting dimensions of $G_P$ are both equal to $1$. Moreover, the projection onto the first coordinate is Lipschitz and maps $G_P$ onto $[0,1]$, so $\dimH G_P\ge1$. Since the Hausdorff dimension does not exceed the upper box-counting dimension (see \cite[Chapters~2 and~3]{Fal14}), we obtain $\dimH G_P=\dimB G_P=1$.

It remains to determine the covering dimension. To apply Lemma~\ref{lem:partitions}, we construct finite partitions of $G_P$ into open sets whose diameters tend to $0$. For $n\in\N$ and $0\le k<2^n$, let $J_{n,k}$ be the union of $I_{n,k}$ and those endpoints $d$ of $I_{n,k}$ for which either $d\in\{0,1\}$ or $P(d)$ is the one-sided limit of $P$ from within $I_{n,k}$. Each point $d\coloneqq j2^{-n}$ with $0<j<2^n$ belongs to $\D$. By Theorem~\ref{thm:main}(\ref{thm:main(i)}), its two one-sided limits are distinct, and $P(d)$ equals exactly one of them. Hence, $d$ lies in exactly one of $J_{n,j-1}$ and $J_{n,j}$. More precisely, we include a point $p/2^m\in\D$ in the interval on its left if $p\equiv1\pmod4$, and in the interval on its right if $p\equiv3\pmod4$. It follows that the intervals $J_{n,k}$, $0\le k<2^n$, form a partition of $[0,1]$. Whether a point $d\in\D$ is included in the interval on its left or on its right depends only on $d$, and not on $n$. The two halves of $I_{n,k}$ are $I_{n+1,2k}$ and $I_{n+1,2k+1}$. If the left endpoint of $I_{n,k}$ belongs to $J_{n,k}$, then it also belongs to $J_{n+1,2k}$, and if the right endpoint of $I_{n,k}$ belongs to $J_{n,k}$, then it also belongs to $J_{n+1,2k+1}$. Moreover, the midpoint $(2k+1)2^{-n-1}$ of $I_{n,k}$ belongs to exactly one of $J_{n+1,2k}$ and $J_{n+1,2k+1}$. Hence, $J_{n,k}$ is the disjoint union of $J_{n+1,2k}$ and $J_{n+1,2k+1}$. Therefore, the sets
\[
G_{n,k}\coloneqq\{(x,P(x)):x\in J_{n,k}\}\qquad(0\le k<2^n)
\]
form a partition $\mathscr P_n$ of $G_P$, and $\mathscr P_{n+1}$ refines $\mathscr P_n$.

We claim that each $G_{n,k}$ is both closed and open in $G_P$. To see that it is closed, let $(y_j,P(y_j))_{j\in\N}$ be a sequence in $G_{n,k}$ converging to a point $(y,P(y))$ of $G_P$. Since $y_j\in J_{n,k}\subseteq[k2^{-n},(k+1)2^{-n}]$ and $y_j\to y$, the point $y$ lies in this closed interval. If $y\notin J_{n,k}$, then $y$ is an endpoint of $I_{n,k}$ other than $0$ and $1$, hence a point of $\D$, and $y_j\in I_{n,k}$ for all sufficiently large $j$. It follows that $P(y_j)$ tends to the one-sided limit of $P$ from within $I_{n,k}$. By the definition of $J_{n,k}$, this limit is different from $P(y)$. This is a contradiction, and so $y\in J_{n,k}$. Thus, $G_{n,k}$ is closed in $G_P$. Now, the other members of $\mathscr P_n$ are also closed in $G_P$, and there are only finitely many of them. Their union is therefore closed in $G_P$. Since $G_{n,k}$ is the complement of this union in $G_P$, it is open in $G_P$, as claimed.

Moreover, by \eqref{eq:cell}, the values of $P$ on $I_{n,k}$ lie in a closed interval of length $(n+2)2^{-n}$. If an endpoint $d$ of $I_{n,k}$ lies in $J_{n,k}$, then $P(d)$ is the one-sided limit from within $I_{n,k}$ by the definition of $J_{n,k}$ and the continuity at $0$ and $1$. Thus, $P(d)$ lies in the same closed interval. It follows that $G_{n,k}$ is contained in a rectangle of width $2^{-n}$ and height $(n+2)2^{-n}$, and its diameter is at most $(n+3)2^{-n}$. Since this bound tends to $0$ as $n\to\infty$, Lemma~\ref{lem:partitions} gives $\dimcov G_P=0$.
\end{proof}

\subsection{Proof of Theorem~\ref{thm:takagi}}\label{subsec:takagi}
\begin{proof}[Proof of Theorem~\ref{thm:takagi}]
The first equality follows from \eqref{eq:remainder}, since $|e_n(x)|=2^{-n}T^n(x)$ for every $x\in[0,1]$ and $n\in\N$. This equality is also valid when the expansion has terminated, since both sides are $0$. By induction on $n$, we have
\[
T^n(x)=2\dist(2^{n-1}x,\Z)\qquad(x\in[0,1],\ n\in\N),
\]
which gives the second equality. The resulting series is the usual definition of the Takagi function (see, e.g., \cite{AK12}).

The series consists of continuous functions and converges uniformly, so $A$ is continuous. We determine the Hausdorff and box-counting dimensions by a covering argument similar to the one in Section~\ref{subsec:iii}. Let $n\in\N$, and consider the sum of the first $n$ terms of the last series in Theorem~\ref{thm:takagi}. Each term is $1$-Lipschitz, so this partial sum is $n$-Lipschitz. Moreover, the sum of the remaining terms lies between $0$ and $2^{-n}$. Hence, for each $0\le k<2^n$, the part of $G_A$ above $[k2^{-n},(k+1)2^{-n}]$ lies in a rectangle of width $2^{-n}$ and height $(n+1)2^{-n}$, and so it can be covered by $n+1$ squares of side $2^{-n}$. Thus, $N_{2^{-n}}(G_A)\le(n+1)2^n$, so the upper box-counting dimension of $G_A$ is at most $1$. As before, the projection onto $[0,1]$ gives $\dimH G_A\ge1$, and hence $\dimH G_A=\dimB G_A=1$. Finally, the map $x\mapsto(x,A(x))$ is a homeomorphism from $[0,1]$ onto $G_A$, so $\dimcov G_A=\dimcov[0,1]=1$ (see \cite{Eng95}).
\end{proof}

\section{Further remarks}\label{sec:remarks}
In this section, we discuss the effect of changing the values at dyadic rationals and give some alternative arguments for the results proved above.

\begin{remark}\label{rem:representative}
If we use terminating binary expansions in \eqref{eq:intro-binary}, then we obtain a right-continuous function $Q$ on $[0,1)$ such that $Q(r)=P(r^+)$ for every $r\in\D$. Put $Q(1)\coloneqq 0$. Then $Q(x)=E_2(x)-x$ for $x\in[0,1)$, and $Q$ differs from $P$ precisely at the points $p/2^m\in\D$ with $p\equiv1\pmod4$, at which $P-Q$ takes the value $2^{1-m}$. Changing the value at a dyadic rational does not remove the jump there.

The graphs $G_P$ and $G_Q$ have the same closure. To see this, note that they agree over $(0,1)\setminus\D$, and each one-sided limit can be obtained by approaching through points of $(0,1)\setminus\D$. The values at $0$ and $1$ are also obtained in this way. Hence, both graphs are contained in the closure of their common part over $(0,1)\setminus\D$. Since a finite union of closed squares is closed, $N_\delta(G_P)=N_\delta(G_Q)$ for every $\delta>0$. Their Hausdorff dimensions also agree, because the graphs differ by countably many points and both project onto $[0,1]$. Finally, we have $Q(r)=Q(r^+)\ne Q(r^-)$ for every $r\in\D$. By taking limits through points of $(0,1)\setminus\D$, we see that \eqref{eq:cell} is also valid for $Q$. Hence, we can apply the proof of $\dimcov G_P=0$ in Section~\ref{subsec:iii} to $Q$. Thus, all assertions of Theorem~\ref{thm:main}(\ref{thm:main(iii)}) hold for $Q$. We point out that the equality of the covering numbers follows from the equality of the closures. For the covering dimension, we also use the fact that the value at each dyadic rational is one of the two one-sided limits. The countability of the set where $P$ and $Q$ differ is not enough to give these conclusions.

We next consider L\'evy's function $\mathcal L$. By Remark~\ref{rem:base2}, it agrees with $P$ on $[0,1]\setminus\D$ and takes the average of the two one-sided limits at each point of $\D$. To see that \eqref{eq:cell} is also valid for $\mathcal L$, fix $I_{n,k}$ and approximate the value at a dyadic point in $I_{n,k}$ by the average of values at points on its two sides. Consequently, the covering estimates in Section~\ref{subsec:iii} also hold for $\mathcal L$. The argument on differentiability in Section~\ref{subsec:ii} applies as well. Indeed, at a point of continuity, it uses only non-dyadic points, while at a point of discontinuity there is no finite derivative. Since $\mathcal L(r)$ lies strictly between its two one-sided limits for every $r\in\D$, the point $(r,\mathcal L(r))$ is isolated in $G_{\mathcal L}$. Hence, $G_{\mathcal L}$ and $G_P$ do not have the same closure. To determine its covering dimension, modify the partitions in Section~\ref{subsec:iii} as follows. For each $n\in\N$, replace $J_{n,k}$ by the union of $I_{n,k}$ and its endpoints in $\{0,1\}$, and add the points of $G_{\mathcal L}$ above $j2^{-n}$, $0<j<2^n$, as singletons. These singletons are open by the preceding observation. The other members of the partition are open by the same argument using one-sided limits as in Section~\ref{subsec:iii}. The partitions still refine one another and satisfy the same bound on the diameters. Hence, Lemma~\ref{lem:partitions} gives $\dimcov G_{\mathcal L}=0$. Therefore, $\mathcal L$ has no finite derivative at any point of $(0,1)$, and all assertions of Theorem~\ref{thm:main}(\ref{thm:main(iii)}) hold for $\mathcal L$. We remark that the assertion on differentiability does not follow from the fact that the pointwise H\"older exponent of $\mathcal L$ is at most $1$ everywhere (see the proof of \cite[Proposition~4]{Jaf97}), since a function may be differentiable at a point where its pointwise H\"older exponent equals $1$.
\end{remark}

\begin{remark}\label{rem:traces}
We can also prove the lower bound for the covering number and the infinite variation by using the one-sided limits. Let $n\in\N$, and write $I_{n,k}=(u,v)$. Since $T^{n+1}(u)=T^{n+1}(v)=0$, every $e_j$ with $j>n$ tends to $0$ as a point approaches either endpoint. On the other hand, $S_n$ has slope $n$ on $I_{n,k}$. By uniform convergence, we therefore obtain $P(v^-)-P(u^+)=n2^{-n}$. By Lemma~\ref{lem:ivp}, applied at points of $I_{n,k}$ sufficiently close to $u$ and $v$, every value strictly between $P(u^+)$ and $P(v^-)$ is attained on $I_{n,k}$. Hence, the vertical projection of the part of $G_P$ above $I_{n,k}$ contains an interval of length $n2^{-n}$. A square of side $2^{-n}$ meets at most two of these open intervals in the horizontal direction. Summing the lengths of the vertical projections, we obtain the sharper bound $N_{2^{-n}}(G_P)\ge n2^{n-1}$. For the variation, fix a non-degenerate interval $[a,b]\subseteq[0,1]$, and let $m\in\N$ satisfy $(b-a)2^{m-2}\ge1$. Since $[a2^m,b2^m]$ is an interval of length $(b-a)2^m$, it contains at least $(b-a)2^{m-1}-1\ge(b-a)2^{m-2}$ odd integers, and so $[a,b]$ contains at least $(b-a)2^{m-2}$ points $p/2^m\in\D$. Each of them has a jump of size $2^{1-m}$ by Theorem~\ref{thm:main}(\ref{thm:main(i)}), so the sum of these jumps is at least $(b-a)/2$. Summing over all such $m$, we see that the sum of the sizes of all jumps of $P$ in $[a,b]$ is infinite. This gives another proof of infinite total variation on $[a,b]$.
\end{remark}

\begin{remark}\label{rem:clopen}
We give another proof of $\dimcov G_P=0$ by using the binary representation and some basic topological facts. Let $K$ be the closure of $G_P$, and equip $\{0,1\}^\N$ with the product topology, where $\{0,1\}$ is discrete. Define $\pi\colon\{0,1\}^\N\to\R^2$ by
\[
\pi(\omega)\coloneqq\biggl(\sum_{j\in\N}\frac{\omega_j}{2^j},\ \sum_{j\in\N}(j-2)\frac{\omega_j}{2^j}\biggr)\qquad\bigl(\omega\coloneqq(\omega_j)_{j\in\N}\in\{0,1\}^\N\bigr).
\]
The two series defining $\pi$ converge uniformly, so $\pi$ is continuous. By \eqref{eq:intro-binary}, if the first coordinate is a non-dyadic point $x$, then the second coordinate is $P(x)$. At a point of $\D$, the two binary expansions give its two one-sided limits, which are distinct by Theorem~\ref{thm:main}(\ref{thm:main(i)}). These are the only possible pairs of distinct sequences with the same first coordinate. Hence, $\pi$ is injective. Since $\{0,1\}^\N$ is compact and $\R^2$ is Hausdorff, $\pi$ is a homeomorphism onto its image.

Sequences which are not eventually constant are dense in $\{0,1\}^\N$, and their images form the graph of $P$ over $(0,1)\setminus\D$. Every point of $G_P$ is in the closure of this graph, by Theorem~\ref{thm:main}(\ref{thm:main(i)}). Therefore, the image of $\pi$ is precisely $K$. The Cantor space $\{0,1\}^\N$ has covering dimension $0$, and so $\dimcov K=0$. Since the covering dimension of a metrizable space is at least that of any of its subspaces (see \cite{Eng95}), we conclude that $\dimcov G_P=0$.

We can also use $\pi$ to see that $K$ is self-affine. In fact, we have
\[
\pi(\omega)=M\pi(\omega_2,\omega_3,\dots)+\omega_1z,\qquad
M\coloneqq\frac12\begin{pmatrix}1&0\\1&1\end{pmatrix},\quad
z\coloneqq\frac12\begin{pmatrix}1\\-1\end{pmatrix}.
\]
Thus, $K=MK\cup(MK+z)$ is a self-affine set (see \cite[Chapter~9]{Fal14}). Here, $MK$ and $MK+z$ are the images under $\pi$ of the sequences with $\omega_1=0$ and $\omega_1=1$, respectively. Since $\pi$ is injective, these two sets are disjoint. For $n\in\N$ and $w\in\{0,1\}^n$, the part of $G_P$ in the image under $\pi$ of the sequences beginning with $w$ is precisely the set $G_{n,k}$ in Section~\ref{subsec:iii}, where $k\coloneqq\sum_{j=1}^nw_j2^{n-j}$. Thus, the partitions used in Section~\ref{subsec:iii} can also be obtained from the binary representation in this way.
\end{remark}

\begin{remark}\label{rem:unsigned}
Recall that Mandelbrot \cite{Man82} defined a fractal as a set whose covering dimension is strictly less than its Hausdorff dimension. By Theorems~\ref{thm:main} and~\ref{thm:takagi}, we have $\dimcov G_P=0<1=\dimH G_P$, while $\dimcov G_A=\dimH G_A=1$. Thus, the graph of the signed error sum is a fractal in the sense of Mandelbrot, whereas the graph of the unsigned error sum is not.
\end{remark}


\begin{thebibliography}{99}
\bibitem{Ahn23} M.~W.~Ahn, {\it On the error-sum function of Pierce expansions}, J. Fractal Geom. {\bf 10} (2023), 389--421.
\bibitem{AK12} P.~C.~Allaart \& K.~Kawamura, {\it The Takagi function: a survey}, Real Anal. Exchange {\bf 37} (2011/12), 1--54.
\bibitem{BBDK96} J.~Barrionuevo, R.~M.~Burton, K.~Dajani, \& C.~Kraaikamp, {\it Ergodic properties of generalized L\"uroth series}, Acta Arith. {\bf 74} (1996), 311--327.
\bibitem{DT11} M.~Dai \& L.~Tang, {\it On the error-sum function of tent map base series}, J. Math. Anal. Appl. {\bf 378} (2011), 571--577.
\bibitem{Eng95} R.~Engelking, {\it Theory of Dimensions, Finite and Infinite}, Sigma Ser. Pure Math. 10, Heldermann, Lemgo, 1995.
\bibitem{Fal14} K.~Falconer, {\it Fractal Geometry: Mathematical Foundations and Applications}, 3rd ed., John Wiley \& Sons, Chichester, 2014.
\bibitem{Jaf97} S.~Jaffard, {\it Old friends revisited: the multifractal nature of some classical functions}, J. Fourier Anal. Appl. {\bf 3} (1997), 1--22.
\bibitem{Jaf04} S.~Jaffard, {\it On Davenport expansions}, in: Fractal Geometry and Applications: A Jubilee of Beno\^it Mandelbrot, Part~1 (M.~L.~Lapidus and M.~van Frankenhuijsen, eds.), Proc. Sympos. Pure Math. 72, Amer. Math. Soc., Providence, RI, 2004, 273--303.
\bibitem{JLZ25} S.~Jaffard, L.~Liao, \& Q.~Zhang, {\it Pointwise regularity and multifractal analysis of pure jump series}, preprint, 2025. Available at \url{https://perso.math.u-pem.fr/jaffard.stephane/pdf/PureJumps.pdf}.
\bibitem{JLZ26a} S.~Jaffard, L.~Liao, \& Q.~Zhang, {\it Multivariate multifractal analysis of L\'evy functions. Part~I: Determination of multifractal spectra}, Proc. Roy. Soc. Edinburgh Sect. A, published online, 2026, doi:10.1017/prm.2025.10120.
\bibitem{JLZ26b} S.~Jaffard, L.~Liao, \& Q.~Zhang, {\it Multivariate multifractal analysis of L\'evy functions. Part~II: validity of the multifractal formalism}, Math. Z. {\bf 312} (2026), Paper No. 68, 21 pp.
\bibitem{KL95} J.-P.~Kahane \& P.-G.~Lemari\'e-Rieusset, {\it Fourier Series and Wavelets}, Studies in the Development of Modern Mathematics 3, Gordon and Breach, Luxembourg, 1995.
\bibitem{Man82} B.~B.~Mandelbrot, {\it The Fractal Geometry of Nature}, W.~H.~Freeman, San Francisco, 1982.
\bibitem{QD10} F.~Qiao \& M.~Dai, {\it On the error-sum function of expansion in base $p$}, Int. J. Nonlinear Sci. {\bf 9} (2010), 330--334.
\bibitem{RP00} J.~N.~Ridley \& G.~Petruska, {\it The error-sum function of continued fractions}, Indag. Math. (N.S.) {\bf 11} (2000), 273--282.
\bibitem{SW07} L.-M.~Shen \& J.~Wu, {\it On the error-sum function of L\"uroth series}, J. Math. Anal. Appl. {\bf 329} (2007), 1440--1445.
\bibitem{Ton07} X.~Tong, {\it On the error-sum function of decimal expansion}, J. Math. (PRC) {\bf 27} (2007), 637--640.
\end{thebibliography}
\end{document}